\title{ ~~\\ Chebyshev's bias for
composite numbers with restricted
prime divisors}
\author{Pieter Moree}
\documentclass[12pt]{article}
\usepackage{amssymb, latexsym, amsfonts}
\def\@ptsize{2}
\newtheorem{Thm}{Theorem}

\newtheorem{Lem}{Lemma}
\newtheorem{Cor}{Corollary}

\newcommand{\qed}{\hfill $\Box$}

\begin{document}
\date{}
\maketitle
{\def\thefootnote{}
\footnote{\noindent P. Moree: KdV Institute, University of Amsterdam,
Plantage Muidergracht 24, 1018 TV Amsterdam, The Netherlands, e-mail:
moree@science.uva.nl}}
{\def\thefootnote{}
}
\begin{abstract}
\noindent Let $\pi(x;d,a)$ denote the number of primes $p\le x$
with $p\equiv a({\rm mod~}d)$.
Chebyshev's bias is the phenomenon that
`more often' $\pi(x;d,n)>\pi(x;d,r)$, than the other way around,
where $n$ is a
quadratic non-residue mod $d$ and $r$ is a quadratic residue mod $d$.
If $\pi(x;d,n)\ge \pi(x;d,r)$ for every $x$ up to some large
number, then one expects that $N(x;d,n)\ge N(x;d,r)$ for
every $x$. Here $N(x;d,a)$ denotes the number of integers $n\le x$
such that every prime divisor $p$ of $n$ satisfies
$p\equiv a({\rm mod~}d)$. In this paper we develop some tools to
deal with this type of problem and apply them to show that, for
example,
$N(x;4,3)\ge N(x;4,1)$ for every $x$.
In the process we express the so called second order
Landau-Ramanujan constant as an infinite series and show
that the same type
of formula holds true for a much larger class
of constants.
\end{abstract}
\section{Introduction}
Especially for small moduli $d$ primes seem to have a preference for
non-quadratic residue classes mod $d$ over quadratic residue
classes mod $d$. This phenomenon is called Chebysev's bias \cite{chebyshev}.
For example, $\pi(x;3,1)$ does not exceed
$\pi(x;3,2)$ for the first time until $x=608981813029$, as was shown
by Bays and Hudson \cite{bays}. On the other hand Littlewood \cite{littlewood}
has shown that the function $\pi(x;3,2)-\pi(x;3,1)$ has
infinitely many sign changes.
Rubinstein and Sarnak \cite{chebbie} quantified some biases under
the assumption of the Generalized Riemann Hypothesis (GRH) and the assumption
that the non-negative imaginary parts of the nontrivial zeros of all
Dirichlet L-functions are linearly independent over the rationals.
Define $\delta_{q,a_1,a_2}$ to be the logarithmic density of the set
of real numbers $x$ such that the inequality $\pi(x;q,a_1)>\pi(x;q,a_2)$
holds, where the logarithmic density of a set $S$ is
$$\lim_{x\rightarrow \infty}{1\over \log x}\int_{[2,x]\cap S}{dt\over t},$$
assuming the limit exists. Under the aforementioned assumptions
and assuming that $(\Bbb Z/q\Bbb Z)^*$ is cyclic,
Rubinstein and Sarnak showed that $\delta_{q,a_1,a_2}$ always exists
and is strictly positive and, moreover,
that $\delta_{q,n,r}>0.5$ if and only if $n$ is a non-square mod $q$
and $r$ is a square mod $q$. They calculated, amongst others,
that $\delta_{4;3,1}=0.9959\cdots$ and $\delta_{3;2,1}=0.9990\cdots$
Thus Chebyshev's bias is not
only an initial interval phenomenon.
The comparison of the behaviour of primes
lying in various arithmetic progressions is the subject of comparative
prime number theory, which was systematically developped in a
series of papers by Knapowski and Tur\'an, cf. \cite{turanbook,turan}.
More recent references are e.g. \cite{martin,chebbie} and
various papers of J. Kaczorowski. One of the older papers, by
Wintner \cite{wintner}, was inspired, interestingly enough, on
the observation (p. 240) that there is `an apparent
parrallelism between certain problems in celestial mechanics on
the one hand and "wobbly" terms of the explicit formula of
Riemann and Von Mangoldt on the other hand". Wintner, who
has written many papers on celestial mechanics then could put
his expertise in that field to good
use.\\
\indent Put $g_{d,a}(n)=0$ if $n$ has no prime divisor $p$ satisfying
$p\equiv a({\rm mod}~d)$ and $g_{d,a}(n)=1$ otherwise (note
that $g_{d,a}(1)=1$). We let $N(x;d,a)=\sum_{n\le x}g_{d,a}(n)$.
The contribution of the small primes to the growth of
$N(x;d,a)$ is much bigger than to $\pi(x;d,a)$ and
hence we might expect
that if $\pi(x;d,a)\ge \pi(x;d,b)$ up to some reasonable
$x$, then actually $N(x;d,a)\ge N(x;d,b)$ for every $x$.
In general, given two non-negative
multiplicative functions $f$ and $g$ such that $f$ and $g$ are equal
to a positive constant $\tau$ in the primes on average and such that
there is a bias towards $f$ in the sense that $\sum_{p\le x}f(p)
\ge \sum_{p\le x}g(p)$ for all $x$ up to some rather large number,
is it true that $\sum_{n\le x}f(n)\ge \sum_{n\le x}g(n)$ for
{\it every} $x$ ? The asymptotic behaviour of the latter type of sums
is well-understood and so proving that these types of results
are true
{\it asymptotically} is usually not difficult. We can for example
invoke the following classical result due to Wirsing \cite{wirsing2}.
\begin{Thm}
\label{eduard} {\rm (Wirsing \cite{wirsing2})}.
Let $f$ be a multiplicative function
satisfying $0\le f(p^r)\le c_1c_2^r,~c_1\ge 1,~1\le c_2<2$,
and
$\sum_{p\le x}f(p)=(\tau+o(1)){x/\log x}$,
where $\tau,c_1$ and
$c_2$ are constants. Then, as $x\rightarrow \infty$,
$$\sum_{n\le x}f(n)\sim {e^{-\gamma\tau}\over \Gamma(\tau)}
{x\over \log x}\prod_{p\le x}\left(1+{f(p)\over p}+
{f(p^2)\over p^2}+{f(p^3)\over p^3}+\cdots+\right),$$
where $\gamma$ is Euler's constant and $\Gamma(\tau)$ denotes the
gamma-function. (Here and in the sequel the letter $p$ is used
to indicate primes.)
\end{Thm}
We thus see that, for $i=1$ and $i=2$,
\begin{equation}
\label{beginnetje}
N(x;3,i)\sim {e^{-\gamma/2}\over \sqrt{\pi}}
{x\over \log x}\prod_{p\le x\atop p\equiv i({\rm mod~}3)}
\left(1-{1\over p}\right)^{-1},
\end{equation}
showing clearly the strong influence of the smaller primes.
By \cite[Theorem 2]{williams} we deduce from the latter formula that
$N(x;3,i)\sim C_{3,i}x/\sqrt{\log x}$, with
\begin{equation}
\label{cdrieeen}
C_{3,1}={3^{1/4}\over \pi\sqrt{2}}\prod_{p\equiv 1({\rm mod~}3)}
\left(1-{1\over p^2}\right)^{-1/2}=
{\sqrt{2}\over 3^{{5\over 4}}}
\prod_{p\equiv 2({\rm mod~}3)}\left(1-{1\over p^2}\right)^{1\over 2},
\end{equation}
where in the derivation of the (\ref{cdrieeen}) we used Euler's
identity $\pi^2/6=\prod_p(1-p^{-2})^{-1}$ (another, selfcontained,
derivation of (\ref{cdrieeen}) is given in Section \ref{constants}).
Using Merten's theorem or \cite[Theorem 2]{williams} again, we
easily infer that $C_{3,2}=2/(3\pi C_{3,1})$.
Restricting to the primes $p\le 29$, we compute that
$C_{3,1}<0.302$ and $C_{3,2}>0.703$ (for more precise
numerical evaluations see Section \ref{constants}).
We thus infer
that $N(x;3,2)\ge N(x;3,1)$ for every sufficiently large $x$.
If we want to make this effective, the extensive literature,
cf. \cite{postnikov},
on multiplicative functions satisfying
conditions as in Wirsing's theorem seems
to offer no help as nobody seems to have been
concerned with proving effective results in this area, which is precisely
what the Chebyshev bias problem for composites challenges us to do.
In this paper we develop some tools for this and apply them to prove:
\begin{Thm}
\label{drietje}
The inequalities $N(x;3,2)\ge N(x;3,1)$, 
$N(x;4,3)\ge N(x;3,1)$,
$N(x;3,2)\ge N(x;4,1)$
and $N(x;4,3)\ge N(x;4,1)$ hold true for every $x$.
\end{Thm}
Not surprisingly the Chebyshev
 bias problem for composites is rather computational in
nature and this appears to preclude one from proving more general
results.\\
\indent The counting functions appearing in Theorem \ref{drietje} can be shown
to satisfy more precise asymptotic estimates than (\ref{beginnetje}).
Theorem \ref{oud} together with the prime number theorem for arithmetic
progressions, shows that there exist constants
$C_{d,a},C_{d,a}(1),C_{d,a}(2),\cdots$ such that for each integer
$m\ge 0$ we have
$$N(x;d,a)={C_{d,a} x\over \sqrt{\log x}}
\left(1+\sum_{j=0}^m{C_{d,a}(j)\over \log^j x}+
O\left({1\over \log^{m+1}x}\right)\right),$$
where the implied constant may depend on $m,a$ and $d$. Thus $N(x;d,a)$
satisfies an asymptotic expansion in the sense of Poincar\'e in terms
of $\log x$. The most famous example of such an expression states that
for $B(x)$, the counting function of the integers that can be represented
as a sum of two integer squares, we have
\begin{equation}
\label{folklore}
B(x)={Kx\over \sqrt{\log x}}
\left(1+\sum_{j=0}^m{K_{j+2}\over \log^j x}+
O\left({1\over \log^{m+1}x}\right)\right),
\end{equation}
where $K$ is the Landau-Ramanujan constant and $K_2$ the second order
Landau-Ramanujan constant.
The Landau-Ramanujan constant is named after Landau
\cite{landau1} who proved in
1908, using
contour integration, that $B(x)\sim Kx/\sqrt{\log x}$ and
Ramanujan, who in his first letter
to Hardy claimed he could prove that
$B(x)=K\int_2^x{dt/\sqrt{\log t}}+O(x^{1/2+\epsilon})$, cf.
\cite{cazaran}. Ramanujan's
claim implies $K_2=1/2$ by partial integration, which was
shown to be false by Shanks \cite{shanks}.
Indeed, we have
$$K={1\over \sqrt{2}}\prod_{p\equiv 3({\rm mod~}4)}
\left(1-{1\over p^2}\right)^{-1/2}=0.76422365358922066299069873125\cdots$$
and
$$K_2={1\over 2}-{\gamma\over 4}-{L'(1,\chi_4)\over 4L(1,\chi_4)}
+{\log 2\over 4}+{1\over 2}
\sum_{p\equiv 3({\rm mod~}4)}{\log p\over p^2-1}=0.5819486593172907\cdots$$
These
constants have been calculated with 1000D precision at least, see
\cite{finch}. It was a folklore result that $B(x)$ should satisfy
(\ref{folklore}), which was written
down by Serre \cite{serre}, who gave some nice
applications to fourier coefficients of modular forms as well.\\
\indent Let $f$ be a non-negative multiplicative function.
Suppose there exists a positive constant $\tau$ such
that
$$\sum_{n\le x}f(n)=\lambda_1(f) x\log^{\tau-1}x\left(1+
(1+o(1)){\lambda_2(f)\over \log x}\right),~x\rightarrow\infty.$$
We then define $\lambda_2(f)$ to be the {\it generalized second-order
Landau-Ramanujan constant}.
In Theorem \ref{tweedeorde} we will identify
a subclass of multiplicative functions for which
this constant exists and
express it as an infinite series.\\
\indent The
second-order generalized Landau-Ramanujan
constant $\lambda_2(f)$ is closely related to the constant $B_f$
appearing in the proof of Lemma \ref{lemma1}
(the key lemma
in the proof of Theorem \ref{drietje}). (I suggest to read the next
section first before reading further.) Lemma \ref{lemma1} yields an
effective estimate for $\mu_f(x)$, provided we can find constants
$\tau,C_{-}$ and $C_{+}$ satisfying (\ref{inklemming}). For the
functions $f$ associated to the quantities in Theorem \ref{drietje}
we find admissible values
of these constants in Section \ref{taudifference}, which requires
effective estimates for counting functions of squarefree numbers of
a certain type (Section \ref{effectivesquarefree}). (At the end
of Section \ref{taudifference} we show that under GRH finding
$C_{-}$ and $C_{+}$ is much easier.) In Section \ref{vier} we show
how to obtain effective estimates for $M_f(x)$ from effective
estimates for $\mu_f(x)$. In Section \ref{kinderachtig} we show how
to prove certain subcases of Theorem \ref{drietje} for every $x$ up to some large $x_0$ using
existing numerical work on the associated Chebyshev prime biases. All
these ingredients then come together in Section \ref{eindelijk}, where
a proof of Theorem \ref{drietje} is given.\\
\indent In Section \ref{vijf} we find an infinite series expansion for
the constant $B_f$ appearing in Lemma \ref{lemma1a} (we have
$C_{-}\le B_f\le C_{+}$) and relate it to the generalized second-order
Landau-Ramanujan constant. Section \ref{constants} contains a numerical study
of some of the constants appearing in this paper.\\
\indent In \cite{hexagonal} the methods developed in this paper are
somewhat refined and then used to resolve
Schmutz Schaller's conjecture (see \cite[p. 201]{schmutz} or
the introduction of \cite{conway}) that the hexagonal lattice
is "better" than the square lattice. More precisely, let
$0<h_1<h_2<\cdots$ be the positive integers, listed in ascending order,
which can be written as $h_i=x^2+3y^2$ for integers  $x$ and $y$. Let
$0<q_1<q_2<\cdots$ be the positive integers, listed in ascending order,
which can be written as $q_i=x^2+y^2$ for integers $x$ and $y$. Then
Schmutz Schaller's conjecture is that $q_i\le h_i$ for $i=1,2,3,\cdots$.

\section{Notation}
Let $f$ be a
non-negative
real-valued multiplicative function.
We define $M_f(x)=\sum_{n\le x}f(n)$,
$\mu_f(x)=\sum_{n\le x}f(n)/n$
and $\lambda_f(x)=\sum_{n\le x}f(n)\log n$.
We denote the formal Dirichlet series $\sum_{n=1}^{\infty}f(n)n^{-s}$
associated to $f$ by $L_f(s)$.
If $f(p)$ equals $\tau>0$ on average at primes
$p$, it can be shown that
$\lim_{s\rightarrow 1+0} (s-1)^{\tau}L_f(s)$ exists, under
some mild additional conditions on $f$.
In that case we
put $$C_f:={1\over \Gamma(\tau)}\lim_{s\rightarrow 1+0} (s-1)^{\tau}L_f(s).$$
We have $C_f>0$. We define $\Lambda_f(n)$
by
$$-{L_f'(s)\over L_f(s)}=\sum_{n=1}^{\infty}{\Lambda_f(n)\over n^s}.$$
Notice that
\begin{equation}
\label{convolutie}
f(n)\log n=\sum_{d|n}f(d)\Lambda_f({n\over d}).
\end{equation}
The notation suggests that $\Lambda_f(n)$ is an analogue of the Von Mangoldt
function. Indeed, if $f={\bf 1}$, then $L_f(s)=\zeta(s)$ and
$\Lambda_f(n)=\Lambda(n)$. From (\ref{convolutie}) we infer by M\"obius
inversion the well-known formula
\begin{equation}
\label{dieisoudzeg}
\Lambda(n)=\sum_{d|n}\mu(d)\log{n\over d}.
\end{equation}
In general, on writing $L_f(s)$ as an Euler product, one
easily sees that $\Lambda_f(n)$ is zero if $n$ is not a prime power.
If $f$ is the
characteristic function of  a subsemigroup
of the natural integers with $(1<)q_1<q_2<\cdots$ as generators,
then it can be shown that
$\Lambda_f(n)=\log q_i$ if $n$ equals a positive power of a generator $q_i$
and $\Lambda_f(n)=0$ otherwise. Thus for example, if
$f=g_{d,a}$, then $\Lambda_{g_{d,a}}(n)=\log p$ if $n=p^r$, $r\ge 1$
and
$p\equiv a({\rm mod~}d)$, and $\Lambda_{g_{d,a}}(n)=0$ otherwise.\\
\indent From property (\ref{convolutie}) of
$\Lambda_f(n)$, we easily infer that
\begin{equation}
\label{lnaarpsi}
\lambda_f(x)=\sum_{n\le x}f(n)\psi_f({x\over n}),
\end{equation}
where $\psi_f(x)=\sum_{n\le x}\Lambda_f(n)$.
For some further properties of $\Lambda_f(n)$ the reader is referred
to \cite[\S 2.2]{cazaran}.\\
\indent The notation $x_0,\alpha$ and $\beta$ is used to indicate
inessential local constants, their values might be different in different
contexts.

\section{Effective estimates for $\mu_f(x)$}
\label{twee}
The following result will play a crucial r\^ole. It
uses some ideas from the proof of Theorem A in \cite{song}.
\begin{Lem}
\label{lemma1}
Let $f$ be a non-negative multiplicative arithmetic function.
Suppose that there exists constants
$\tau(>0),C_{-}$ and $C_+$
such that
\begin{equation}
\label{inklemming}
C_{-}\le \sum_{n \le x}{\Lambda_f(n)\over n}-\tau \log
x \le C_{+}~~{\rm for~every~}x\ge 1,
\end{equation}
then, for $x>\exp(C_{+})$, we have
\begin{equation}
\label{firstestimate}
{C_f\over \tau}\log^{\tau}x{\left(1-{C_+\over \log x}\right)^{\tau+1}\over
1-{C_{-}\over \log x}}
\le \mu_f(x)
\le {C_f\over \tau}\log^{\tau}x
{\left(1-{C_{-}\over \log x}\right)^{\tau+1}\over
1-{C_{+}\over \log x}}
\end{equation}
where
\begin{equation}
\label{defcee}
C_f:={1\over \Gamma(\tau)}\lim_{s\rightarrow 1+0} (s-1)^{\tau}L_f(s).
\end{equation}
\end{Lem}
Remark. An alternative expression for $C_f$ is given by
$$C_f={1\over \Gamma(\tau)}\lim_{s\rightarrow 1+0}\prod_p\left(
1+\sum_{k=1}^{\infty}{f(p^k)\over p^{ks}}\right)
\left(1-{1\over p^s}\right)^{\tau}.$$
{\it Proof of Lemma} \ref{lemma1}. Let $B_f$ be an arbitrary constant and write
\begin{equation}
\label{startingpoint}
\sum_{n\le x}{\Lambda_f(n)\over n}=\tau \log
x + B_f+E_f(x).
\end{equation}
(This is unnecessary for this proof, but needed in the proof of
Lemma \ref{lemma2}, so we do this now to save some space later.)
We have
\begin{eqnarray}
\sum_{n \le x}{f(n)\log  n \over n}&=&\sum_{ n \le x}\sum_{d|n}{f(d)\over d}{\Lambda_f(n/d)\over {n\over d}}
= \sum_{d \le x}{f(d)\over d}\sum_{k\le {x\over d}}
{\Lambda_f(k)\over k}\nonumber\\
&=&\tau \sum_{n\le x}{f(n)\over n}\log({x\over n})
+ B_f\mu_f(x)+\sum_{ n \le x}{f(n)\over n}E_f({x\over n}).\nonumber
\end{eqnarray}
We write this equality in the form
$$-\sum_{n\le x}{f(n)\over n}\log({x\over n})+\mu_f(x)\log x   =
\tau \sum_{n\le x}{f(n)\over n}\log({x\over n})
+B_f\mu_f(x)+\sum_{n \le x}{f(n)\over n}E_f({x\over n}).$$ This inequality on its turn can
be written, using that
$$\sum_{n\le x}{f(n)\over n}\log{x\over n}
=\sum_{n\le x}{f(n)\over n}\int_n^x{dt\over t}=
\int_1^x{\mu_f(t)\over t}dt,$$
 as
\begin{equation}
\label{functionalequation}
\mu_f(x)\log x - (\tau+1)\int_1^x{\mu_f(v)\over v} dv = B_f\mu_f(x)
+\sum_{n\le x}{f(n)\over n}E_f({x\over
n}).
\end{equation}
Put $\sigma_f(x)=\int_1^x{\mu_f(v)\over v}dv$.
By assumption $C_{-}\le B_f+E_f(x)\le C_{+}$ for $x\ge 1$.
Using
(\ref{functionalequation}) we then deduce that
$\mu_f(x)=(\tau+1)\sigma_f(x)/\log x+\mu_f(x)\epsilon_f(x)$,
where $C_{-}\le \epsilon_f(x)\log x\le C_+$.
Solving this for $\mu_f(x)$ we find that
\begin{equation}
\label{useful}
\mu_f(x)={1\over 1-\epsilon_f(x)}{\tau+1\over \log x}\sigma_f(x),~~x\ge x_0,
\end{equation}
where $x_0:=\exp((1+\delta)C_+)$, and $\delta>0$ is arbitrary and fixed.
In the rest of the proof we assume that
$x\ge x_0$.
Let
$$R_f(t):=\log\left({\tau+1\over \log ^{\tau+1}t}\sigma_f(t)\right).$$
Note that, for $t\ge x_0$,
\begin{equation}
\label{afgeleideschatting}
R'_f(t)={\tau+1\over t\log t}
{\epsilon_f(t)\over [1-\epsilon_f(t)]}
\end{equation}
and hence $R'_f(t)=O(t^{-1}\log^{-2}t)$.
Thus $\int_x^{\infty}R'_f(t)dt$ converges absolutely, and
therefore $\int_x^{\infty}R'_f(t)dt=A_f-R_f(x)$, for some constant $A_f$
not depending on $x$.
On writing $D_f=\exp(A_f)$ we obtain
\begin{equation}
\label{tweemaalis}
{\tau+1\over \log^{\tau+1}x}\sigma_f(x)=\exp(R_f(x))=
D_f\exp\left(-\int_x^{\infty}R_f'(t)dt\right)
\end{equation}
Using (\ref{afgeleideschatting}) and $C_{-}\le \epsilon_f(x)\log x\le C_+$, we
see that
$$\int_x^{\infty}{C_{-}(\tau+1)\over t\log t[\log t-C_{-}]}dt
\le \int_x^{\infty}R'_f(t)dt\le \int_x^{\infty}
{C_+(\tau+1)\over t\log t[\log t-C_+]}dt.$$
Thus
$$-(\tau+1)\log\left(1-{C_{-}\over \log x}\right)
\le \int_x^{\infty}R'_f(t)dt \le -(\tau+1)\log\left(1-{C_{+}\over
\log x}\right).$$
On combining (\ref{tweemaalis}) with (\ref{afgeleideschatting}) we
deduce that
\begin{equation}
\label{poehpoeh}
D_f \left(1-{C_+\over \log x}\right)^{\tau+1}
\le {\tau+1\over
\log^{1+\tau}x}
\sigma_f(x)\le D_f\left(1-{C_{-}\over \log x}\right)^{\tau+1}
\end{equation}
We will now show that $D_f=C_f/\tau$. The inequalities
(\ref{poehpoeh}) in combination with (\ref{useful}) imply in particular that
\begin{equation}
\label{hulp}
\mu_f(x)=D_f\log ^{\tau} x + O(\log^{\tau-1}x).
\end{equation}
 By partial integration and using the well-known integral expression for
the gamma function we find that
\begin{eqnarray}
L_f(s)&=&(s-1)\int_1^{\infty}{\mu_f(t)\over t^s}dt=(s-1)\int_1^{\infty}{D_f\log ^{\tau} t+O(\log ^{\tau-1} t)\over t^s}dt\nonumber\\
&=&D_f{\Gamma(\tau+1)\over (s-1)^{\tau}}+O\left({s-1\over (s-1)^{\tau}}\right)
\nonumber
\end{eqnarray}
and thus $D_f=C_f/\tau$.
The inequalities (\ref{poehpoeh}) together with (\ref{useful}) yield (\ref{firstestimate})
on using that $D_f=C_f/\tau$ and $C_{-}\le \epsilon_f(x)\log x\le C_+$.
\qed\\

\noindent The convolutional nature of $\sum_{n\le x}E(x/n)f(n)/n$ forces us
to require that $x\ge 1$ in (\ref{inklemming}) (whereas we would like
to replace it with $x\ge x_0$). Nevertheless we can invoke the following
easy lemma to improve on (\ref{firstestimate}).
\begin{Lem}
\label{lemma1a}
Suppose that there exists constants $D_{-}$ and
$D_{+}$ such that for every $x\ge x_0$,
\begin{equation}
\label{volgendeinklemming}
D_{-}\mu_f(x)\le B_f\mu_f(x)+\sum_{n\le x}{f(n)\over n}E_f({x\over n})\le
D_{+}\mu_f(x).
\end{equation}
Then we have, for $x>\max\{x_0,\exp(D_{+})\}$,
\begin{equation}
\label{secondestimate}
{C_f\over \tau}\log^{\tau}x{\left(1-{D_+\over \log x}\right)^{\tau+1}\over
1-{D_{-}\over \log x}}
\le \mu_f(x)
\le {C_f\over \tau}\log^{\tau}x
{\left(1-{D_{-}\over \log x}\right)^{\tau+1}\over
1-{D_{+}\over \log x}}.
\end{equation}
\end{Lem}
{\it Proof}. Follows easily on closer scrutiny of the
previous proof. \qed\\

\noindent We now give an example of how Lemma \ref{lemma1a} can be used. By
assumption we have $E_f(x)\le C_{+}-B_f$ for every $x\ge 1$. Suppose
that $E_f(x)\le C_{+}'-B_f$ for $x\ge n_0$, where $C_{+}'<C_{+}$.
An upper bound for the innerterm in (\ref{volgendeinklemming}) is then given by
$$C_{+}\mu_f(x)-(C_{+}-C'_{+})\mu_f({x\over n_0}).$$
Using the explicit bounds in
(\ref{firstestimate})
we can then find an $x_0$ and
$D_{+}<C_{+}$ such that the conditions of Lemma \ref{lemma1a}
are satisfied (note
that $D_{+}>C'_{+}$). By applying (\ref{secondestimate}) instead
of (\ref{firstestimate}) a better value for $D_{+}$ can then
be obtained. Then iterate.\\

\noindent By making an
assumption on $E_f(x)$ we will,
not surprisingly, be able to do better
than both Lemma \ref{lemma1} and Lemma \ref{lemma1a}.
\begin{Lem}
\label{lemma2}
Let $f$ be a non-negative multiplicative arithmetic function and suppose
that
{\rm (\ref{startingpoint})} holds with
\begin{equation}
\label{egbound}
|E_f(x)|\le {c_0\over {\rm max}\{1,\log x\}}.
\end{equation}
for every $x\ge 1$, where $c_0$ is some explicit constant. Then
there exist effectively computable
constants $\alpha,\beta$ and $x_0$ such that
$$\mu_f(x)={C_f\over \tau}\log^{\tau}x-C_fB_f\log^{\tau-1}x+
{\cal E}_f(x),$$
where $\alpha\log^{\tau-1/2}x\le {\cal E}_f(x)
\le \beta\log^{\tau-1/2}x$ for every $x\ge x_0$.
\end{Lem}
{\it Proof}. We denote the right hand side
of (\ref{egbound}) by $h(x)$ and put
$s(x)=x/e^{\sqrt{\log x}}$. Let $x_0\ge e$.
Using Lemma \ref{lemma1} it is not difficult to see that
$$\beta_f(x_0):=\sup_{x\ge x_0}\sqrt{\log x}\left\{1-
{\mu_f(s(x))\over \mu_f(x)}
\left(1-{1\over \sqrt{\log x}}\right)\right\},$$
is finite and can be effectively computed (note
that $\beta_f(x_0)\ge \tau+1$).\\
\indent Clearly
$$\Big|\sum_{n\le x}{f(n)\over n}E_f({x\over n})\Big|
\le \sum_{n\le s(x)}{f(n)\over n}h({x\over n})
+\sum_{s(x)<n\le x}{f(n)\over n}h({x\over n}).$$
Denote the latter two sums by $I_1$ and $I_2$.
We have $I_1\le c_0\mu_f(s(x))/\sqrt{\log x}$ and
$I_2\le c_0(\mu_f(x)-\mu_f(s(x))$.
We thus find that (\ref{useful}) holds true with
$$|\epsilon_f(x)-{B_f\over \log x}|
\le {c_0\beta_f(x_0)\over \log^{3/2} x},~x\ge x_0.$$
Proceeding as in
the proof of Lemma \ref{lemma1}, but with this improved error
estimate, the result then easily follows. \qed\\

\section{Relating $\mu_f(x)$ to $M_f(x)$}
\label{vier}
\indent Given an effective
estimate for $\mu_f(t)$, we can derive an
effective estimate for $M_f(t)$ on using that
\begin{equation}
\label{omhoog}
M_f(x)-M_f(x_0)=\int_{x_0}^x t~d\mu_f(t).
\end{equation}
Suppose that
$${C_f\over \tau}\log^{\tau}x \left(1+{\alpha\over \log x}\right)
\le \mu_f(x) \le
{C_f\over \tau}\log^{\tau}x \left(1+{\alpha+\beta\over \log x}\right),$$
for some constants $\alpha$ and $\beta$ and every $x\ge x_0$ (if
the conditions of Lemma \ref{lemma1} are satisfied, such
$\alpha,\beta$ and $x_0$ can certainly be determined). This
leads to an upperbound for $M_f(x)$ that is asymptotically equal
to $C_f(1+\beta/\tau)x\log^{\tau-1}x$ and a lowerbound that is
asymptotically equal to $\max\{0,C_f(1-\beta/\tau)x\log^{\tau-1}x\}$. These
estimates are too weak for our purposes.\\
\indent Write
$\mu_f(x)=C_f\log^{\tau}x/\tau-C_fB_f\log^{\tau-1}x+
{\cal E}_f(x),$ cf. Lemma \ref{lemma2}, and suppose that
${\cal E}^{-}_f(x)\le {\cal E}_f(x)\le {\cal E}^{+}_f(x)$
for every $x\ge x_0$, where ${\cal E}^+_f(x)$ and ${\cal E}^{-}_f(x)$
are effectively computable.
(This supposition is
certainly true if the conditions of Lemma \ref{lemma2} are
satisfied.)
Let $C_f(x_0)=M_f(x_0)-x_0{\cal E}_f(x_0)-{C_f}x_0\log^{\tau}x_0$.
Then
an easy computation shows that for every $x\ge x_0$,
$$M_f(x)\le C_fx\log^{\tau-1}x+(1-\tau)C_f(1+B_f)
\int_{x_0}^x{\log^{\tau-2}t}~dt+C_f(x_0)+R_f(x),$$
where
$$
x{\cal E}^{-}_f(x)-\int_{x_0}^x{\cal E}^{+}_f(t)~dt\le R_f(t)
\le x{\cal E}^{+}_f(x)-\int_{x_0}^x{\cal E}^{-}_f(t)~dt.$$
There are various problems with this approach, one of the major
ones being getting a good estimate for $c_0$ in Lemma \ref{lemma2}.\\
\indent An alternative approach starts with the observation that, for
$x\ge 2$,
\begin{equation}
\label{mlambda}
M_f(x)=\int_{2-}^x{d\lambda_f(t)\over \log t}=
{\lambda_f(x)\over \log x}+\int_2^x
{\lambda_f(t)\over t\log ^2 t}~dt,
\end{equation}
and that if we have explicit bounds of the type $\alpha x<\psi_f(x)<\beta x$,
then $\lambda_f(x)$ can be related to $x\mu_f(x)$ by (\ref{lnaarpsi}). Note
in particular that if $\lambda_f(x)\ge \lambda_g(x)$ for every $x\ge 2$, 
then $M_f(x)\ge M_g(x)$ for every $x$ (the reverse implication is not
always true in general). The disadvantage of proving something stronger
is hopefully compensated by the fact that $\lambda_f(x)$ can be
easily related to $\mu_f(x)$.

\section{The generalized second-order Ramanujan-\\ Landau constant}
\label{vijf}
In Theorem \ref{tweedeorde} we will identify
a subclass of multiplicative functions for which
the generalized Landau-Ramanujan constant
(defined in Section 1) exists and relate it to
an infinite series involving $\Lambda_f(n)$.
The following result will play an essential r\^ole in this.
\begin{Thm}
\label{oud}
{\rm \cite[Theorem 6]{cazaran}}.
Let $f$ be a multiplicative function
satisfying
\begin{equation}
\label{grens}
0\le f(p^r)\le c_1c_2^r,~c_1\ge 1,~1\le c_2<2,
\end{equation}
and
$\sum_{p\le x}f(p)=\tau {\rm Li}(x)+O\left({x\log^{-2-\rho}x}\right),$
where $\tau$ and $\rho$ are positive real fixed numbers.
Then there exists a constant $B_f$ such that
{\rm (\ref{startingpoint})} holds with
$E_f(x)=O(\log^{-\rho}x)$.
Moreover, for
every $\epsilon>0$,
\begin{equation}
\label{ouwetje!}
\sum_{n\le x}{f(n)\over n}=\sum_{0\le \nu <\rho+1}a_{\nu}
\log^{\tau-\nu}x+O(\log^{\tau-1-\rho+\epsilon} x),
\end{equation}
where the implied constant depends at most on $f$ and $\epsilon$.
In case $f$ is completely multiplicative, condition {\rm (\ref{grens})} can
be weakened to
$\sum_{p,r\ge 2,~p^r>x}{(f(p)/p)^r\log p}=O(\log^{-\rho}x)$.
\end{Thm}
{\it Proof}. This result is just Theorem 6 of \cite{cazaran}, except for the
claim regarding $E_f(x)$, the truth of which is
however established in the course of the proof of Theorem 6
of \cite{cazaran}. \qed\\

\noindent The next result shows that the second-order Landau-Ramanujan constant
is closely related to the constant $B_f$ appearing in (\ref{startingpoint}).
\begin{Thm}
\label{tweedeorde}
Let $f$ be a multiplicative function
satisfying the hypothesis of Theorem {\rm \ref{oud}}
with $\rho>1$.
Then $\lambda_2(f)$, the generalized second-order Landau-Ramanujan constant,
equals
$$\lambda_2(f)=(1-\tau)\left(1+\tau\gamma+\sum_{n=1}^{\infty}{\Lambda_f(n)-\tau\over n}
\right),$$
or alternatively
$\lambda_2(f)=(1-\tau)(1+B_f)$,
where
$$B_f:=\lim_{x\rightarrow \infty}\left(\sum_{n\le x}{\Lambda_f(n)\over n}
-\tau \log x\right).$$
\end{Thm}
{\it Proof}.
Since by assumption $\rho>1$, we have by (\ref{ouwetje!})
$$\mu_f(x)=a_0\log^{\tau} x +a_1 \log^{\tau-1}x
+a_2\log^{\tau-2}x+O(\log^{\tau-2-\delta}x),$$
for some $\delta>0$.
Theorem \ref{oud} implies that
$B_f$
exists. Using that $\log x=\sum_{n\le x}1/n-\gamma+o(1),$ we
see that it suffices to prove that $\lambda_2(f)=(1-\tau)(1+B_f)$.
Theorem
\ref{oud} yields that $E_f(x)=O(\log^{-1}x)$, hence the conditions
of Lemma \ref{lemma2} are satisfied and it follows
that $a_0=\tau C_f$ and $a_1=-C_f B_f$.
On using that $M_f(x)=x\mu_f(x)-
\int_1^x\mu_f(t)dt$ it follows by partial integration
that $\lambda_1(f)=C_f$ and $\lambda_2(f)=
(1-\tau)(1+B_f)$,
as required. \qed\\

\noindent Example. Let $b_1$ be the characteristic function of
the set of natural numbers that can be written as a sum of two integer
squares. This is a subsemigroup of the natural numbers that is generated
by the primes $p$ with $p\equiv 1({\rm mod~}4)$, $p=2$ and the
squares of the remaining prime numbers (this result goes back to Fermat).
By what has been said in Section 2, it then follows that
$$
\Lambda_{b_1}(n)=\cases{
2\log p &if $n=p^{r},~r\ge 1$ and $p\equiv 3({\rm mod~}4)$;\cr
\log p &if $n=p^r,~r\ge 1$ and $p\equiv 1({\rm mod~}4)$ or $p=2$;\cr
0 &otherwise.}
$$
Application of Theorem \ref{tweedeorde} yields the following two
formulae for
the second-order Landau-Ramanujan constant $K_2$ (cf. Section 1):
$$K_2={1\over 2}\left(1+{\gamma\over 2}+\sum_{n=1}^{\infty}
{\Lambda_{b_1}(n)-{1\over 2}\over n}\right)=
{1\over 2}\lim_{x\rightarrow \infty}\left(
1+\sum_{n\le x}{\Lambda_{b_1}(n)\over n}-{1\over 2}\log x\right).$$

\section{Numerical evaluation of certain constants}
\label{constants}
In order to complete our proof we need to evaluate certain
constants with enough precision. For some of them this has
been done before.\\
\indent We first consider the evaluation of $C_{3,1}$ and
$C_{3,2}$. We have, for $\Re (s)>1$,
$L_{g_{3,1}}(s)=\prod_{p\equiv 1({\rm mod~}3)}
(1-p^{-s})^{-1}$. Note that
\begin{equation}
\label{generating}
L_{g_{3,1}}(s)^{2}=\zeta(s)L(s,\chi_3)(1-3^{-s})
\prod_{p\equiv 2({\rm mod~}3)}
(1-p^{-2s}).
\end{equation}
From this, (\ref{defcee}), $\lim_{s\rightarrow 1+0}(s-1)\zeta(s)=1$ and
the fact that $\Gamma({1\over 2})=\sqrt{\pi}$, we obtain
$$C_{3,1}^2={2L(1,\chi_3)\over 3\pi}
\prod_{p\equiv 2({\rm mod~}3)}
\left(1-{1\over p^2}\right).$$
If $\chi$ is a real primitive character modulo $k$ and $\chi(-1)=-1$,
then
$$L(1,\chi)=-{\pi\over k^{3/2}}\sum_{n=1}^k n\chi(n),$$
by Dirichlet's celebrated class number
 formula (cf.
equation (17) of Chapter 6 of \cite
{davenport}).
We infer that $L(1,\chi_3)=\pi/\sqrt{27}$.
Using that $C_{g_{3,1}}\ge 0$ and $\zeta(2)=\pi^2/6$, we
then deduce (\ref{cdrieeen}). Using that
$L_{g_{3,2}}(s)L_{g_{3,1}}(s)(1-3^{-s})^{-1}=\zeta(s)$, we infer
that $C_{3,2}=2/(3\pi C_{3,1})$.\\
\indent In order to compute $C_{3,2}$ and $C_{3,1}$ with many decimal
accuracy we proceed as in Shanks \cite[p. 78]{shanks}.
We note that, for $\Re (s)>1/2$,
\begin{equation}
\label{shanksje}
\prod_{p\equiv 2({\rm mod~}3)}(1-p^{-2s})^{2}
={L(2s,\chi_3)\over \zeta(2s)(1-3^{-2s})}\prod_{p\equiv 2({\rm mod~}3)}
(1-p^{-4s}),
\end{equation}
from which we infer by recursion that
$$C_{3,1}={\sqrt{2}\over 3^{5\over 4}}\prod_{n=1}^{\infty}\left(
{L(2^n,\chi_3)\over (1-3^{-2^n})
\zeta(2^n)}\right)^{1\over 2^{n+1}}.$$
Because of the lacunary character of this expression, it can be
calculated quickly up to high precision, which yields
$C_{3,1}=0.3012165544749342124\cdots$ and
$C_{3,2}=0.7044984335\cdots$.
Similarly one can show that $C_{4,3}=1/(2\pi C_{4,1})$ and
$$C_{4,1}={1\over 2\sqrt{2}}\prod_{p\equiv 3({\rm mod~}4)}
\left(1-{1\over p^2}\right)^{1/2}=
{1\over \pi}\prod_{p\equiv 1({\rm mod~}4)}
\left(1-{1\over p^2}\right)^{-1/2}.$$
Using Shanks' trick we then infer
that $C_{4,1}=0.3271293669410263824002328\cdots$
and $C_{4,3}=0.4865198883\cdots$.\\
\indent On noting that, for $\Re (s)\ge 1$,
$$\sum_{n=1}^{\infty}{\Lambda(n)-1\over n^s}
=-{\zeta'(s)\over \zeta(s)}-\zeta(s),$$
and using that $\zeta(s)=1/(s-1)+\gamma+O(s-1)$ is the Taylor series for $\zeta(s)$ around
$s=1$ (see e.g. \cite[pp. 162-164]{narkiewicz}), one infers that
\begin{equation}
\label{lambdaovern}
\sum_{n\le x}{\Lambda(n)\over n}=\sum_{n\le x}{1\over n}
-2\gamma+o(1)=\log x -\gamma+o(1).
\end{equation}
Taking the logarithmic derivative of
(\ref{generating}) one obtains that
$$-2{L_{g_{3,1}}'\over L_{g_{3,1}}}(s)=
-{\zeta'\over \zeta}(s)-{L'\over L}(s,\chi_3)
-{\log 3\over 3^s-1}-2\sum_{p\equiv 2({\rm mod~}3)}{\log p\over
p^{2s}-1},$$
from which one easily infers that
$$2\sum_{n\le x}{\Lambda_{g_{3,1}}(n)\over n}
=\sum_{n\le x}{\Lambda(n)\over n}-{L'\over L}(1,\chi_3)
-{\log 3\over 2}-2\sum_{p\equiv 2({\rm mod~}3)}{\log p\over
p^{2}-1}+o(1),$$
which yields, on invoking (\ref{lambdaovern}),
$$2B_{g_{3,1}}=-\gamma-{L'\over L}(1,\chi_3)-{\log 3\over 2}
-2\sum_{p\equiv 2({\rm mod~}3)}{\log p\over p^2-1}.$$
Similarly we deduce that
$$2B_{g_{4,1}}=-\gamma-{L'\over L}(1,\chi_4)
-{\log 2}-2\sum_{p\equiv 3({\rm mod~}4)}{\log p\over p^2-1}.$$
As to the numerical evaluation of $B_{g_{3,1}}$ and
$B_{g_{4,1}}$, we note that
$$2\sum_{p\equiv 2({\rm mod~}3)}{\log p\over p^2-1}=-
{d\over ds}\log \prod_{p\equiv 2({\rm mod~}3)}\left({1\over 1-p^{-2s}}
\right)\Big|_{s=1}.$$
Then, applying (\ref{shanksje}) $m$ times, we obtain
$$
\sum_{p\equiv 2({\rm mod~}3)}{\log p\over p^2-1}=
\sum_{p\equiv 2({\rm mod~}3)}{\log p\over p^{2^{m+1}}-1}
+{1\over 2}\sum_{n=1}^m\left\{{L'\over L}(2^m,\chi_3)
-{\zeta'\over \zeta}(2^m)-{\log 3\over 3^{2^m}-1}\right\}.$$
Now $L$-functions and their derivatives can be computed with
high accuracy using for example PARI
(cf. \cite[Section 10.3]{cohenbook}).
On doing so we find that the prime sum in the left
hand side of the latter formula equals
$0.3516478132638087560157790\cdots$.
Similarly we have
$$
\sum_{p\equiv 3({\rm mod~}4)}{\log p\over p^2-1}=
\sum_{p\equiv 3({\rm mod~}4)}{\log p\over p^{2^{m+1}}-1}
+{1\over 2}\sum_{n=1}^m\left\{{L'\over L}(2^m,\chi_4)
-{\zeta'\over \zeta}(2^m)-{\log 2\over 2^{2^m}-1}\right\}.$$
We thus find that the sum on the left hand side equals
$0.2287363531940324576\cdots$.
For more on evaluating infinite sums or products involving primes,
we refer to \cite{cohenpreprint} and \cite{moree}.\\
\indent For the
logarithmic derivative
 $L'(1,\chi)/L(1,\chi)$ we find, with $\chi=\chi_3$
 and $\chi=\chi_4$,
${L'\over L}(1,\chi_3)=
0.36828161597014784263323790407578664254876430999
\cdots$ and
$${L'\over L}(1,\chi_4)=
0.2456095847773141723888166261790625184335337829549\cdots$$
An alternative way of evaluating the latter two
logarithmic derivatives is by relating them to the gamma
function or the arithmetic-geometric-mean (AGM).
We have (Berger (1883), Lerch (1897), de S\'eguier (1899)
and Landau \cite{landau0}),
$${L'\over L}(1,\chi_4)=
\log\left(M(1,\sqrt{2})^2{e^{\gamma}\over 2}\right),$$
where $M(1,\sqrt{2})$ denotes the
limiting value of Lagrange's AGM algorithm
$a_{n+1}=(a_n+b_n)/2$, $b_{n+1}=\sqrt{a_nb_n}$ with inputs $a_0=1$ and
$b_0=\sqrt{2}$.
It can be shown that $M(1,\sqrt{2})=\sqrt{2\over \pi}
\Gamma({3\over 4})^2$.
Gauss showed (in his diary), cf. \cite{cox}, that
$${1\over M(1,\sqrt{2})}={2\over \pi}\int_0^1{dx\over \sqrt{1-x^4}}.$$
The total arclength of the lemniscate $r^2=\cos(2\theta)$ is given
by $2L$, where
$L:=\pi/M(1,\sqrt{2})$ is the
so-called {\it lemniscate constant}. If $\chi=\chi_3$ we have
similarly, with $z:=\sin({\pi\over 12})={(\sqrt{3}-1)\over \sqrt{8}}$,
$${L'\over L}(1,\chi_3)=\log
\left({2^{4\over 3}M(1+z,1-z)^2 e^{\gamma}\over 3}\right),~
M(1+z,1-z)={2^{4\over 3}\pi^2\over 3^{1\over 4}\Gamma({1\over 3})^3},$$
and
$${1\over M(1+z,1-z)}={3^{1\over 4}\over \pi}\int_0^1
{dx\over \sqrt{x(1-x^3)}}.$$
The values of $L'(1,\chi_4)$ and $L'(1,\chi_3)$ can also be determined
using generalized Euler constants for arithmetical progressions, see
Examples 1 and 2 of \cite{dilcher}.
For general non-trivial real $\chi$ the  quotients $L'(1,\chi)/L(1,\chi)$
`feel' the zeroes of $L(s,\chi)$ close to 1
(see \cite[pp. 80-83]{davenport} for a quantitative version)
and a study of
their average behaviour might throw some light on the (non)-existence
of the Landau-Siegel zeros, cf. \cite{blackholes}.\\
\indent
Putting
our
subcomputations together, we find
that
$$B_{g_{3,1}}=-1.09904952586667653048446536830561\cdots$$
On noting that  $\Lambda(n)=\Lambda_{g_{3,1}}(n)
+\Lambda_{g_{3,2}}(n)$ if for $3\nmid n$, it is easily deduced that
$B_{g_{3,1}}=-\gamma-{\log 3\over 2}-B_{g_{3,2}},$
using which we compute that
$$B_{g_{3,2}}=-.02747228336891117581966934023805\cdots$$
Similarly, we find
$B_{g_{4,1}}=-0.9867225683134286288516284\cdots$
and
$$B_{g_{4,3}}=-\log 2-\gamma-B_{g_{4,1}}=
-0.2836402771480495411721157\cdots$$

\indent An alternative approach in calculating the constant
 $B_{g_{d,a}}$ is on invoking the formula
$$\varphi(d)\sum_{n\le x\atop n\equiv a({\rm mod~}d)}
{\Lambda(n)\over n}=\log x-\gamma-\sum_{p|d}{\log p\over p-1}-
\sum_{\chi\ne \chi_0}\chi({\bar a}){L'\over L}(1,\chi)+o(1),$$
where $a$ and $d$ are coprime integers, the sum is over the characters
mod $d$ different from the principal character and $\bar a$ is any
integer such that $a{\bar a}\equiv 1({\rm mod~}d)$.
The latter formula is
derived by elementary means in \cite{nevanlinna}.\\
\indent Using Theorem \ref{tweedeorde} we are now in the position
to compute some second-order Landau-Ramanujan constants. They are simply
given by $\lambda_2(f)=(1+B_f)/2$ for
$f\in \{g_{3,1},g_{3,2},g_{4,1},g_{4,3}\}$.\\

\section{Effective estimates for squarefree integers}
\label{effectivesquarefree}
\noindent In the
sequel we will establish some effective estimates
for certain number theoretic functions of a real variable. The
general procedure is to establish the estimates for every $x\ge
x_0$ for some $x_0$. The following lemma can then often be used
to show that there exists
a number $x_1<x_0$ such
that the estimates in fact hold true for every
$x\ge x_1$. It reduces a seemingly continuous problem to a discrete one.
\begin{Lem}
\label{contnaareindig}
Let $y_1>y_0$ be arbitrary
real numbers. Let $F$ and $r$ be non-decreasing
real-valued functions such that, moreover, $F$ changes
its value only at integers. Let $x_1,x_2,\cdots,x_n$ be the integers
in $(y_0,y_1)$ where $F$ changes its value.
Put $x_0=y_0$ and $x_{n+1}=y_n$.
Then
$$\sup_{y_0\le x\le y_1}\{F(x)-r(x)\}=\max_{0\le i\le n}\{F(x_i)-r(x_i)\}$$
and $$\inf_{y_0\le x\le y_1}\{F(x)-r(x)\}=
\min_{0\le i\le n}\{F(x_i)-r(x_{i+1})\}.$$
\end{Lem}
In our proof of Theorem \ref{drietje} we need effective estimates for
$Q_{\chi_3}(x)$ and $Q_{\chi_4}(x)$, where
$Q_{\chi}(x)$ denotes the number of integers $n\le x$ such that
$\mu(n)\chi(n)\ne 0$.
Note that $Q_{\chi_4}$ merely counts the odd
squarefree numbers and hence we will use the
more suggestive notation $Q_{odd}$ for it.
There are two
obvious approaches in estimating these
functions; relating them to $Q(x)$, where $Q(x)$ denotes the number
of squarefree integers not exceeding $x$, and an ab initio approach. We
demonstrate both approaches.\\
\indent Put
$R(x)=Q(x)-6x/\pi^2$. It was shown by Moser and MacLeod
\cite{moser} that
$|R(x)|<\sqrt{x}$ for all $x$ and that
$|R(x)|<\sqrt{x}/2$ for $x\ge 8$.
Cohen and Dress \cite{cohen} showed that $|R(x)|<0.1333\sqrt{x}$ for $x\ge 1664$.
\begin{Lem}
\label{uitoudeschattingen}
For $x\ge 0$ we have
$$\left|Q_{\chi_3}(x)-{9\over 2\pi^2}x\right|\le 0.3154\sqrt{x}
+17.2,$$
and
\begin{equation}
\label{gemakkelijk}
\left|Q_{\chi_3}(x)-{9\over 2\pi^2}x\right|\le {1\over 2}\sqrt{x}+1.
\end{equation}
\end{Lem}
{\it Proof}.
We clearly have $Q(x)=Q_{\chi_3}(x)+Q_{\chi_3}(x/3)$, from
which we infer that
\begin{equation}
\label{drieid}
Q_{\chi_3}(x)=\sum_{i=0}^{\infty}(-1)^iQ({x\over 3^i}).
\end{equation}
Put $x_0=1664$. On applying Lemma \ref{contnaareindig} with
$y_0=x_0/3$ and $y_1=x_0$, we find
that $|R(x)|\le 0.15\sqrt{x}$ in the
interval $(x_0/3,x_0]$. Similarly
we compute that $|R(x)|\le 0.29\sqrt{x}$ in
the interval $(x_0/27,x_0/3]$. These estimates yield
when combined with identity (\ref{drieid}) and the quoted
bounds for $|R(x)|$:
$$\left| Q_{\chi_3}(x)-{9\over 2\pi^2}x\right|\le
\alpha{\sqrt{3x}\over \sqrt{3}-1}+(0.15-\alpha)\sqrt{x_0}+
(0.29-\alpha)\left(\sqrt{{x_0\over 3}}
+\sqrt{{x_0\over 9}}\right)$$
$$+
(0.5-\alpha)\left(\sqrt{{x_0\over 27}}+\sqrt{x_0\over 81}\right)
+(1-\alpha)(\sqrt{9}+\sqrt{3}+1+{1\over \sqrt{3}}+{1\over 3}+
\cdots),$$
where $\alpha=0.1333$. The latter bound does not exceed
$0.3154\sqrt{x}+17.2$.
From this bound we then infer that (\ref{gemakkelijk})
holds  for every $x\ge 10000$.
We now apply
Lemma \ref{contnaareindig} with $y_0=0$ and $y_1=10000$ to
establish the validity of  (\ref{gemakkelijk})
 in the remaining range.\qed\\
\indent Using that
that $|R(x)|\le 0.15\sqrt{x}$ in the
interval $(x_0/2,x_0]$ and $|R(x)|\le 0.29\sqrt{x}$ in
the interval $(x_0/32,x_0/2]$ we deduce, proceeding as
in the proof of Lemma \ref{uitoudeschattingen}, that
$|Q_{odd}(x)-4x/\pi^2|\le 0.4552\sqrt{x}+26.5$. Although the latter
bound is sharp enough for our purposes, we present a selfcontained
proof of a slightly sharper bound (which uses ideas from
\cite{moser}).
\begin{Lem}
\label{qodd}
For $x \ge 0$ we have
$$\left|Q_{odd}(x)-{4\over \pi^2}x\right|
\le {1\over 2}\sqrt{x}+1{\rm ~and~}\left|Q_{odd}(x)-{4\over \pi^2}x\right|
\le ({2\over \pi^2}+{1\over 4})\sqrt{x}+{1\over 4}x^{1\over 4}
+2.$$
\end{Lem}
{\it Proof}. We have
$$Q_{odd}(x)=\sum_{n\le x\atop n{\rm ~odd}}|\mu(n)|=
\sum_{n\le x\atop n{\rm ~odd}}\sum_{d^2|n}\mu(d)
=\sum_{d\le x\atop d{\rm ~odd}}\mu(d)\left[{x\over 2d^2}+{1\over 2}\right].$$
Put $R_{odd}(x)=Q_{odd}(x)-4x/\pi^2$.
On noting that $\sum_{d{\rm ~odd}}\mu(d)/d^2=8/\pi^2$, we find that
$$
|R_{odd}(x)|\le \left|\sum_{d^2\le x\atop d{\rm ~odd}}
\mu(d)\left({x\over 2d^2}-\left[{x\over 2d^2}+{1\over 2}\right]
\right)\right|+x\left|\sum_{d^2>x\atop d{\rm ~odd}}{\mu(d)\over 2d^2}\right|.
$$
Since $|x-[x+1/2]|\le 1/2$ for every $x$, we deduce that
\begin{equation}
\label{laternognodig}
|R_{odd}(x)|\le
{Q_{odd}(\sqrt{x})\over 2}+x\left|\sum_{d^2>x\atop d{\rm ~odd}}{\mu(d)\over 2d^2}\right|.
\end{equation}
Suppose that $x>4$, then
$$\sum_{d^2>x\atop d{\rm ~odd}}{1\over d^2}\le
\sum_{m>{\sqrt{x}-1\over 2}}{1\over (2m+1)(2m-1)}
=\sum_{m>{\sqrt{x}-1\over 2}}\left[{1\over 4m-2}-{1\over 4m+2}\right]
\le {.5\over \sqrt{x}-2}.$$
On using this and the trivial estimate $Q_{odd}(x)\le (x+1)/2$,
we deduce that $|R_{odd}(x)|\le {1\over 2}\sqrt{x}+1$
on applying Lemma \ref{contnaareindig}
with $y_0=0$ and $y_1=36$. Using the latter
bound for $Q_{odd}(x)$ in (\ref{laternognodig})
one then easily obtains the second stated bound in the
formulation of the lemma on applying Lemma \ref{contnaareindig}
with $y_0=0$ and $y_1=9$.

\section{On the difference $\sum_{n\le x}{\Lambda_f(n)\over n}-\tau\log x$}
\label{taudifference}
In order to use Lemma \ref{lemma1} we need to find finite constants $C_{+}$ and
$C_{-}$ such that
$$C_{-}\le \sum_{n\le x}{\Lambda_f(n)\over n}-\tau \log x\le C_{+}$$
for every $x\ge 1$. Recall that $\psi_f(x):=\sum_{n\le x}\Lambda_f(n)$.
Suppose that $\psi_f(x)=\tau x+{\cal  E}_f(x),$ where
$|{\cal  E}_f(x)|\le c_{\epsilon}\log^{-1-\epsilon}x$ for $x\ge x_0$.
Then
$$\sum_{n\le x}{\Lambda_f(n)\over n}=\tau \log x+B_f
+{{\cal  E}_f(x)\over x}-\int_x^{\infty}{{\cal  E}_f(t)\over t^2}dt,$$
and thus, for $x\ge x_0$,
\begin{equation}
\label{difference}
\left|
\sum_{n\le x}{\Lambda_f(n)\over n}-\tau \log x-B_f\right|
\le {c_{\epsilon}\over \log^{\epsilon}x}\left({1\over \epsilon}
+{1\over \log x}\right).
\end{equation}
For example, if $f={\bf 1}$, it is known that $|\theta(x)-x|\le
3.965x/\log^2 x$ for $x>1$ \cite[p.14]{dusart}. Using this with the
bound $\psi(x)-\theta(x)<1.43\sqrt{x}$ \cite[Theorem 13]{rossers}, we can compute $C_{+}$ and $C_{-}$ in this case.
Instead
of carrying this out along
these lines, we proceed slightly differently as this will result
in a sharper bound for the difference in (\ref{difference}).
\begin{Lem}
\label{rosserachtig}
For $x\ge 97$ we have
$$-{1\over 2\log x}+{1\over 2\sqrt{x}}\le
\sum_{n\le x}{\Lambda(n)\over n}-\log x+\gamma\le {2\over \sqrt{x}}
+{1\over 2\log x}.$$
The upper bound holds even true for every $x>1$.
\end{Lem}
{\it Proof}. By \cite[Theorem 6]{rossers} we have, for $x\ge 319$,
$$\left|\sum_{p\le x}{\log p\over p}-\log x- E\right|<{1\over 2\log x},$$
where $E=-\gamma-\sum_{p}\log p\sum_{k\ge 2}p^{-k}$.
Notice that
$$\sum_{n\le x}{\Lambda(n)\over n}=\sum_{p\le x}{\log p\over p}
+\sum_{k\ge 2}\sum_{p}{\log p\over p^k}-\sum_{p^k>x\atop k\ge 2}{\log p\over
p^k}.$$
By partial integration we find that
\begin{equation}
\label{kwadrisum}
\sum_{p^k>x\atop k\ge 2}{\log p\over
p^k}={\theta(x)-\psi(x)\over x}+\int_x^{\infty}{\psi(t)-\theta(t)\over
t^2}dt.
\end{equation}
Suppose that $\alpha \sqrt{t}\le \psi(t)-\theta(t)\le
\beta \sqrt{t}$ for $t\ge x_0$. Then, for $x\ge x_0$ the sum in
(\ref{kwadrisum}) is in the interval $({2\alpha-\beta\over \sqrt{x}},
{2\beta-\alpha\over \sqrt{x}})$. By Theorems 13 and 14 of
\cite{rossers} we
can take $\alpha=0.98$ and $\beta=1.4262$ when
$x_0=319$. On combining
the various estimates, the result follows after
some numerical analysis in the interval $(1,319)$.\qed\\

\noindent From Lemma \ref{rosserachtig} and Lemma
\ref{contnaareindig}
with $y_0=1$ and $y_1=215$ it is
easily deduced that
$$\sup_{x\ge 1}\left\{\sum_{n\le x}{\Lambda(n)\over n}-\log x\right\}
=-{\log 2\over 2}= -0.34657359\cdots$$
and
$$\inf_{x\ge 1}
\left\{\sum_{n\le x}{\Lambda(n)\over n}-\log x\right\}={\log 2\over 2}
-\log 3= -0.75203869\cdots$$
Other than for $f={\bf 1}$, the author is unaware of cases where
an unconditional effective upper bound for ${\cal  E}_f(x)$ of order
$\log^{-1-\epsilon}x$ is known. Thus in order to obtain
admissible values for $C_{+}$
and $C_{-}$ in the case $f
\in \{g_{3,1},g_{3,2},g_{4,1},g_{4,3}\}$ we have to
follow another approach, which is what we will do now.
Notice that
\begin{equation}
\label{ookdatnog}
2\sum_{n\le x}{\Lambda_{g_{3,1}}(n)\over n}=
\sum_{n\le x}{(1+\chi_3(n))\Lambda(n)\over n}
-2\sum_{p^r\le \sqrt{x}\atop p\equiv 2({\rm mod~}3)}{\log p\over p^{2r}}
-\sum_{1<3^r\le x}{\log 3\over 3^r}.
\end{equation}
The latter two sums are easily explicitly estimated and we already
explicitly estimated $\sum_{n\le x}\Lambda(n)/n$. If we can
explicitly estimate $\sum_{n\le x}\chi_3(n)\Lambda(n)/n$, we are
done then. In order to achieve the latter goal, we need a few lemmas.
\begin{Lem}
\label{moebiusje}
Let $h$ be a completely multiplicative function with $h(1)=1$, then
if $g(x)=\sum_{n\le x}h(n)f({x\over n})$
for every $x$, it follows that
$f(x)=\sum_{n\le x}h(n)\mu(n)g({x\over n})$.
\end{Lem}
{\it Proof}. Substitute the expression $\sum_{mn\le x}h(m)f(x/mn)$ for
$g(x/n)$ in the sum $\sum_{n\le x}h(n)\mu(n)g(x/n)$. The resulting expression
simplifies to $f(x)$.\qed\\
\begin{Lem}
\label{voordrieenvier}
Let $\chi$ be a non-principal character and $m_0>1$ be the smallest
integer $>1$ such that $\chi(m_0)\ne 0$. Then
$$
\sum_{n\le x}{\chi(n)\Lambda(n)\over n}+
{L'\over L}(1,\chi)=
O\left({L'\over L}(1,\chi){m_0\over x}Q_{\chi}({x\over m_0})\right)
+O\left({1\over x}\sum_{d\le x/m_0\atop \mu(d)\chi(d)\ne 0}\log{x\over d}\right).
$$
\end{Lem}
{\it Proof}. On using (\ref{dieisoudzeg}) and writing $n=dd_1$ we obtain, for
an arbitrary character $\chi$,
\begin{equation}
\label{basisje}
\sum_{n\le x}{\chi(n)\Lambda(n)\over n}=
\sum_{d\le x/m_0}{\chi(d)\mu(d)\over d}\sum_{d_1\le x/d}
{\chi(d_1)\log d_1\over d_1}.
\end{equation}
On inserting
$$\sum_{d_1\le x/d}{\chi(d_1)\log d_1\over d_1}=-L'(1,\chi)
+O\left({\log(x/d)\over x/d}\right)$$
in this, we obtain
\begin{equation}
\label{vorige}
\sum_{n\le x}{\chi(n)\Lambda(n)\over n}=-L'(1,\chi)
\sum_{d\le x/m_0}{\chi(d)\mu(d)\over d}
+O\left({1\over x}\sum_{d\le x/m_0\atop \mu(d)\chi(d)\ne 0}\log{x\over d}\right).
\end{equation}
We apply Lemma
\ref{moebiusje} with $h(n)={\chi(n)\over n}$ and $f(n)=1$
together with
$\sum_{n\le x}\chi(n)/n=L(1,\chi)+O(1/x)$ to
the latter equation and obtain
\begin{eqnarray}
1&=&\sum_{n\le x/m_0}{\chi(n)\mu(n)\over n}
\left(L(1,\chi)+O({nm_0\over x})\right)\nonumber\\
&=&L(1,\chi)\sum_{n\le x/m_0}{\chi(n)\mu(n)\over n}
+O\left({m_0\over x}Q_{\chi}({x\over m_0})\right).\nonumber
\end{eqnarray}
Combining the latter equation with (\ref{vorige}) and using
the well-known fact that $L(1,\chi)\ne 0$, the result
then follows. \qed\\

\noindent Remark. By using more refined elementary methods
\cite{nevanlinna} one can show that
actually, as $x$ tends to infinity,
$$\sum_{n\le x}{\chi(n)\Lambda(n)\over n}+{L'\over L}(1,\chi)
=o(1).$$
\noindent Let us consider the case
where $\chi=\chi_3$ or $\chi=\chi_4$.
Then, for $x>0$,
\begin{equation}
\label{eleen}
\left|\sum_{n\le x}{\chi(n)\over n}-L(1,\chi)\right|\le
{1\over x},
\end{equation}
where we use that the
non-zero terms in the sum are alternating in sign
and monotonically decreasing.
The function $\log x/x$ is
only decreasing for $x>e$ and a similar argument
then shows that, for $x>e$,
\begin{equation}
\label{elafgeleide}
\left|\sum_{n\le x}{\chi(n)\log n\over n}+L'(1,\chi)\right|
\le {\log x\over x}.
\end{equation}
A numerical analysis shows, however,
that (\ref{elafgeleide}) is still valid for every $x\ge 2$.
The implication of these estimates is that for these characters
and $x\ge 1$ all the implied constants in the latter lemma and its
proof are $\le 1$.
Note that for $x\ge 1$
$$\sum_{d\le x\atop \chi(d)\mu(d)\ne 0}\log{x\over d}
=\sum_{d\le x\atop \chi(d)\mu(d)\ne 0}\int_d^x{dt\over t}
=\int_1^x{Q_{\chi}(t)\over t}dt$$
and thus, for $x\ge m_0$,
$$\sum_{d\le x/m_0\atop \chi(d)\mu(d)\ne 0}\log{x\over d}
=\int_1^{{x\over m_0}}{Q_{\chi}(t)\over t}dt+Q_{\chi}({x\over m_0})
\log m_0.$$
We thus find that, for $x\ge m_0$,
$$x\left|\sum_{n\le x}{\chi(n)\Lambda(n)\over n}+
{L'(1,\chi)\over L(1,\chi)}\right|\le
\left({L'\over L}(1,\chi)m_0+\log m_0\right)
Q_{\chi}({x\over m_0})+\int_1^{x\over m_0}{Q_{\chi}(t)\over
t}dt.$$
For $\chi=\chi_3$ we see, using (\ref{gemakkelijk}), that the right hand side
is bounded by
$$\left({L'\over L}(1,\chi_3)+{\log 2\over 2}+{1\over 2}\right)
{9\over 2\pi^2}+{1\over \sqrt{2x}}\left(
{L'\over L}(1,\chi_3)
+{\log 2\over 2}+1\right)$$
$$+{\log(x/2)\over x}+{1\over x}\left(2{L'\over L}(1,\chi_3)
+\log 2\right).$$
For $\chi=\chi_4$ we see, using that
$Q_{\chi_4}(t)\le 4t/\pi^2+\sqrt{t}/2+1$ (Lemma
\ref{qodd}), that the right hand side
is bounded by
$$\left({L'\over L}(1,\chi_4)+{\log 3\over 3}+{1\over 3}\right)
{4\over \pi^2}+{1\over \sqrt{3x}}\left(
{3\over 2}{L'\over L}(1,\chi_4)
+{\log 3\over 2}+1\right)$$
$$+{\log(x/3)\over x}+{1\over x}\left(3{L'\over L}(1,\chi_4)
+\log 3\right).$$
In the
case where $\chi=\chi_3$ it remains to 
explicitly estimate the latter two sums in (\ref{ookdatnog}).
We have
$$
\sum_{p^r>\sqrt{x}\atop p\equiv a({\rm mod~}d)}{\log p\over p^{2r}}
\le \sum_{p^r>\sqrt{x}}{\log p\over p^{2r}}
=-{\psi(\sqrt{x})\over x}+2\int_{\sqrt{x}}^{\infty}{\psi(t)\over t^3}dt.$$
Using that $0.8t\le \psi(t)\le 1.04t$ for $t\ge 17$ (this
easily follows from Theorem 10 and Theorem 12
from \cite{rossers}), we find that
\begin{equation}
\label{afschattie}
\sum_{p^r>\sqrt{x}\atop p\equiv a({\rm mod~}d)}{\log p\over p^{2r}}\le {1.3\over \sqrt{x}} {\rm ~for~}x\ge 289.
\end{equation}
Furthermore, for every fixed $v>1$ and every $x>0$,
$${\log v\over v-1}
(1-{v\over x})\le \sum_{1<v^r\le x}{\log v\over v^r}\le {\log v\over v-1},$$
where the sum is over the integral powers of $v$ not exceeding $x$. (These
two estimates can also be used in the case where $\chi=\chi_4$.)\\
\indent Let us define
$$C_{+}(f)=\sup_{x\ge 1}\left(\sum_{n\le x}{\Lambda_f(n)\over n}-
{\tau \log x}\right)=B_f+\sup_{x\ge 1}E_f(x),$$
and let $C_{-}(f)$ be similarly defined, with sup replaced by inf.
Let $\epsilon>0$ be fixed. Note that the sharpest
result the method we followed here allows us to prove, with
enough numerical computation, is 
$$\lim_{x\rightarrow \infty}|E_f(x)|\le \left({L'\over L}(1,\chi_3)+{\log 2\over 2}+{1\over 2}\right)
{9\over 4\pi^2}+\epsilon=0.2769537767\cdots+\epsilon{\rm ~and~}$$
\begin{equation}
\label{beperking}
\lim_{x\rightarrow\infty}
|E_g(x)|\le \left({L'\over L}(1,\chi_4)+{\log 3\over 3}+{1\over 3}\right)
{2\over \pi^2}+\epsilon=0.1915268284\cdots+\epsilon,
\end{equation}
where $f\in \{g_{3,1},g_{3,2}\}$ and $g\in \{g_{4,1},g_{4,3}\}$.\\
\indent On putting the various effective bounds together we arrive
at the following result, after numerical
calculations not going
beyond the interval $[1,10^5]$.
\begin{Thm}
\label{zonderGRH} We have\\
{\rm a)} $C_{-}(g_{4,1})>-1.202$ and $C_{+}(g_{4,1})=0$.\\
{\rm b)} $C_{-}(g_{4,3})={\log 3\over 3}-{\log 7\over 2}
=-0.606750\cdots$
 and $C_{+}(g_{4,3})=0$.\\
{\rm c)} $C_{-}(g_{3,1})>-1.4$ and $C_{+}(g_{3,1})=0$.\\
{\rm d)} $C_{-}(g_{3,2})=-{\log 2\over 2}
=-0.34657\cdots$ and $C_{+}(g_{3,2})<0.2764.$
\end{Thm}
\indent On GRH it is much easier to find
the $C_{+}$ and $C_{-}$
satisfying (\ref{inklemming}), which is what will be demonstrated now.
By RH$(d)$ we indicate the hypothesis that for every character
$\chi$ mod $d$ every non-trivial zero of $L(s,\chi)$ is on the
critical line.
Put $$H(x;d,a):
=\sum_{1<p^r\le x\atop p\equiv a({\rm mod~}d)}{\log p\over p^r}
{\rm ~and~}\psi(x;d,a):=\sum_{n\le x\atop n\equiv a({\rm mod~}d)}\Lambda(n).$$
\begin{Lem}
\label{undergrh}
For $d\le 432$ and $(a,d)=1$,
there exists a constant $c_{d,a}$
such that for $x\ge 224$ we have, on RH$(d)$, that
\begin{equation}
\label{steptwo}
\left|\sum_{n\le x\atop n\equiv a({\rm mod~}d)}
{\Lambda(n)\over n}-{\log x\over \varphi(d)}-c_{d,a}\right|\le {11\over 32\pi\sqrt{x}}
\{3\log^2 x+8\log x+16\},
\end{equation}
\end{Lem}
{\it Proof}. In \cite{dusart} it is proved that
for $d\le 432$ and $x\ge 224$ we have, on RH$(d)$, that
\begin{equation}
|\psi(x;d,a)-{x\over \varphi(d)}|\le {11\over 32\pi}\sqrt{x}\log^{2}x.
\end{equation}
Using the latter estimate and partial integration, the lemma
then follows. \qed\\
\indent Using the latter lemma we can compute, under GRH, the exact
values of $C_{-}(g_{4,1}),C_{+}(g_{3,1})$ and $C_{+}(g_{3,2})$.
\begin{Thm}
\label{cplusafschat} 
We have\\
{\rm a)} $C_{-}(g_{4,1})=H(197;4,1)-{\log(229)\over 2}=-0.99076124051235\cdots$, on RH$(4)$.\\
{\rm b)} $C_{-}(g_{3,1})=H(3121;3,1)-{\log(3163)\over 2}=-1.100304022673\cdots$, on
RH$(3)$.\\
{\rm c)} $C_{+}(g_{3,2})=H(5;3,2)-{\log 5\over 2}=
{3\over 4}\log 2-{3\over 10}\log 5=0.03702\cdots$, on RH$(3)$.
\end{Thm}
{\it Proof}. \\ a) Note
that $c_{4,1}=B_{g_{4,1}}$.  On applying Lemma \ref{undergrh} with $d=4$ and $a=1$, Lemma \ref{contnaareindig}, (\ref{afschattie}) and using the numerical value
for $B_{g_{4,1}}$ given in Section \ref{constants}, we deduce 
that $C_{-}(g_{4,1})=\min_{q_i\le 1.79*10^{9}}(H(v_i;4,1)-\log(v_{i+1})/2),$
where $5=v_1<v_2<\cdots$ are the consecutive prime powers 
$p^r$ with $p\equiv 1({\rm mod~}4)$.\\
b) In this case we have
$C_{-}(g_{3,1})=\min_{q_i\le 2.935*10^{10}}(H(q_i;3,1)-\log(q_{i+1})/2),$
where $7=q_1<q_2<\cdots$ are the consecutive prime powers $p^r$ with
$p\equiv 1({\rm mod~}3)$.\\
c) Now $C_{+}(g_{3,2})=\max_{w_i\le 1582079}(H(w_i;3,2)-\log(w_{i})/2)$, where $2=w_1<w_2<\cdots$ are the consecutive prime powers $p^r$ with
$p\equiv 2({\rm mod~}3)$.\\

\section{Connections with Chebyshev's bias for primes}
\label{kinderachtig}
In this section we make some observations that allow us to prove, for
example, 
that $N(x;3,2)\ge N(x;3,1)$ for every $x\le x_0$ for some
large $x_0$, using known
numerical observations regarding $\pi(x;3,2)$ and $\pi(x;3,1)$.\\
\indent Let $Q_1=\{q_1,q_2,q_3,\cdots\}$  and $Q_2=\{v_1,v_2,v_3,\cdots\}$ be sets of pairwise coprime prime powers that satisfy $q_1<q_2<q_3<\cdots$
and $v_1<v_2<v_3<\cdots$.
Let $S_1$ denote the set of integers 
of the form $q_1^{e_1}\cdots q_s^{e_s}$ with
$q_i\in Q_1$ and $e_i\in \Bbb Z_{\ge 0}$ for $1\le i\le s$. Let $S_2$
be similarly defined, but with $Q_1$ replaced by $Q_2$.
Let $\pi_1(x),~\pi_2(x)$, count the number of elements in $Q_1$,
respectively $Q_2$, up to $x$.
If $n=q_1^{e_1}\cdots q_s^{e_s}\in S_1$, then
$m:=v_1^{e_1}\cdots v_s^{e_s}$ is said to be its {\it associate}
in $S_1$. Let $h:\Bbb N\rightarrow \Bbb R_{\ge 0}$ be a non-increasing
function. 
Put $V_1(x)=\sum_{n\in S_1}h(n)$ and $V_2(x)=\sum_{n\in S_2}h(n)$. In
the rest of this section $x_0$ denotes some arbitrary number.
\begin{Lem}
\label{flauw} We have\\
{\rm a)} If $\pi_1(x)\ge \pi_2(x)$ for $x\ge 0$, then
$V_1(x)\ge V_2(x)$ for $x\ge 0$.\\
{\rm b)}  If $\pi_1(x)\ge \pi_2(x)$ for $x\le x_0$, then
$V_1(x)\ge V_2(x)$ for $x\le x_0$.
\end{Lem}
{\it Proof}. a) The assumption implies that if $m\in S_2$, then its associate
$n\in S_1$ satisfies $n\le m$ and $h(n)\ge h(m)$. Thus clearly $V_1(x)\ge V_2(x)$.
The proof of part b will be obvious to the reader now. \qed
\begin{Cor}
If $\pi(x;d,a)\ge \pi(x;d,b)$ for $x\le x_0$, then for $x\le x_0$ we
have both $N(x;d,a)\ge N(x;d,b)$ and $\mu_{g_{d,a}}(x)
\ge \mu_{g_{d,b}}(x)$.
\end{Cor}
The hypothesis in the corollary is in general not strong enough to infer
that $\lambda_{g_{d,a}}(x)\ge \lambda_{g_{d,b}}(x)$ if $x\le x_0$.
However, we have the following easy result.
\begin{Lem}
\label{alsjemenou}
If $M_f(x)\ge M_g(x)$ and $\psi_f(x)\ge \psi_g(x)$ for every $x\le x_0$,
then $\lambda_f(x)\ge \lambda_g(x)$ for $x\le x_0$.
\end{Lem}
{\it Proof}. Use (\ref{lnaarpsi}). \qed
\begin{Cor}
If $$\pi(x;d,a)\ge \pi(x;d,b){\rm ~and~}\sum_{1<p^r\le x\atop p\equiv
a({\rm mod~}d)}\log p\ge \sum_{1<p^r\le x\atop p\equiv
b({\rm mod~}d)}\log p$$ for
every $x\le x_0$, then  
$\lambda_{g_{d,a}}(x)\ge \lambda_{g_{d,b}}(x)$ for $x\le x_0$.
\end{Cor}
In the proof of Theorem \ref{sterker} we will put Corollary 2 a few
times to action.

\section{The proof of Theorem \ref{drietje}}
\label{eindelijk}
The proof of Theorem \ref{drietje} will easily follow from the
following theorem.
\begin{Thm} 
\label{sterker}
For every $x$ we have $\lambda_{g_{3,2}}(x)\ge \lambda_{g_{3,1}}(x)$,
 $\lambda_{g_{3,2}}(x)\ge \lambda_{g_{4,1}}(x)$ and
 $\lambda_{g_{4,3}}(x)\ge \lambda_{g_{3,1}}(x)$.
For  $x\ge 7$ we have $\lambda_{g_{4,3}}(x)\ge \lambda_{g_{4,1}}(x)$.
\end{Thm}
Note that
$$e^{\lambda_{g_{d,a}}(x)}=\prod_{n\le x\atop p|n\Rightarrow p\equiv a({\rm mod~}d)}n.$$
In the proof of Theorem \ref{sterker} we will make use of the following
lemma.
\begin{Lem}
\label{psiafschat} We have
$\psi_{g_{3,1}}(x)\le 0.50456x$ for $x\ge 0$,
$\psi_{g_{3,2}}(x)\ge 0.335x$ for $x\ge 5$,
$\psi_{g_{4,1}}(x)\le 0.50456x$ for $x\ge 0$ and
$\psi_{g_{4,3}}(x)\ge 0.48508x$ for $x\ge 127$.
\end{Lem}
{\it Proof}. Let $d\le 13$ and $(a,d)=1$. Then 
$|\psi(x;d,a)-x/\varphi(d)|\le \sqrt{x}$ 
for $224\le x\le 10^{10}$
by \cite[Theorem 1]{ramare}
and $|\psi(x;d,a)-{x\over \varphi(d)}|<0.004560{x\over \varphi(d)}$ for $x\ge 10^{10}$
by \cite[Theorem 5.2.1]{ramare}.
From these inequalities the lemma 
follows after some computation. \qed \\

\noindent In our proof we consider inequalities of the form
\begin{equation}
\label{veelvariabelen}
\log^{\tau}({x\over r}){\left(1-{C_{+}\over \log(x/r)}\right)^{\tau+1}\over
\left(1-{C_{-}\over \log(x/r)}\right)}
\ge c_1 \log^{\tau}({x\over s}){\left(1-{C'_{-}\over \log(x/s)}\right)^{\tau+1}\over
\left(1-{C'_{+}\over \log(x/s)}\right)},
\end{equation}
where all variables and constants are real numbers with $\tau,r,s$ and
$c_1$ positive, $C_{-}\le C_{+}$, $C'_{-}\le C'_{+}$ and
$x\ge x_0:=\max\{\exp(C'_{+})s,\exp(C_{+})r\}$.
This inequality can be rewritten as
\begin{equation}
\label{rewritten}
1+{C'_{-}-C_{+}+\log(s/r)\over \log(x/s)-C'_{-}}
\ge \left[c_1{\left(1+{C'_{+}-C'_{-}\over \log(x/s)-C'_{+}}\right)
\over {1+{C_{-}-C_{+}\over \log(x/r)-C_{-}}}
}\right]^{1\over \tau}.
\end{equation}
Note that for $x>x_0$ the right hand side is a non-increasing function
of $x$. If $C'_{-}+\log s\le C_{+}+\log r$, the left hand side is
non-decreasing, whereas if the latter inequality is not satisfied the
left hand side asymptotically decreases to 1. We thus arrive at
the following conclusion.
\begin{Lem}
\label{delaatstehoopik}
If $\log s+C'_{-}\le C_{+}+\log r$ and {\rm (\ref{veelvariabelen})} is satisfied
for some $x_1>x_0$, then {\rm (\ref{veelvariabelen})} is satisfied for every
$x\ge x_1$. If $\log s+C'_{-}>C_{+}+\log r$, and the right hand side
of {\rm (\ref{rewritten})} does not exceed 1 for some $x_1>x_0$, then
{\rm (\ref{veelvariabelen})} is satisfied for every $x\ge x_1$.
\end{Lem} 
\noindent {\it Proof of Theorem} \ref{sterker}.\\
$\lambda_{g_{3,2}}(x)$ {\tt versus} $\lambda_{g_{3,1}}(x)$.
Using Lemma \ref{psiafschat} we infer that
$$\lambda_{g_{3,2}}(x)\ge \sum_{n\le {x\over 5}}g_{3,2}(n)\psi_{g_{3,2}}({x\over n})\ge 0.335 \mu_{g_{3,2}}({x\over 5}), {\rm ~and~that}$$
$$\lambda_{g_{3,1}}(x)=\sum_{n\le x}g_{3,1}(n)\psi_{g_{3,1}}({x\over n})=\sum_{n\le {x\over 7}}g_{3,1}(n)\psi_{g_{3,1}}({x\over n})\le 0.50456\mu_{g_{3,1}}({x\over 7}).$$
With $d=3$, $a=2$ and $b=1$ the conditions of Corollary 2 are satisfied for every $x<196699$ (but not for
$x=196699$ as $\psi_{g_{3,1}}(196699)>\psi_{g_{3,2}}(196699)$).
Thus we certainly may assume that $x>1900$.
Using the estimates $C_{3,1}<0.302$ and $C_{3,2}>0.703$ we then deduce
from Lemma 1, Theorem \ref{zonderGRH} and Lemma \ref{delaatstehoopik} that
$0.335\mu_{g_{3,2}}(x/7)>0.50456\mu_{g_{3,1}}(x/7)$.\\
$\lambda_{g_{3,2}}(x)$ {\tt versus} $\lambda_{g_{4,1}}(x)$.
The conditions of Corollary 2 are now satisfied for every $x\le 10^7$
(the smallest $x$ for which the conditions are not satisfied is not known, but
must be less than $10^{82}$ by \cite{fordhudson}). Thus we certainly may assume that $x>4600$.
Then reasoning as before we infer that
$\lambda_{g_{3,2}}(x)\ge 0.335\mu_{g_{3,2}}(x/5)\ge 0.50456\mu_{g_{4,1}}(x/5)\ge \lambda_{g_{4,1}}(x)$.\\
$\lambda_{g_{3,2}}(x)$ {\tt versus} $\lambda_{g_{4,1}}(x)$. 
The conditions of Corollary 2 are now satisfied for every $x\le 10^7$
(the smallest $x$ for which the conditions are not satisfied is not known, but
must be less than $10^{82}$ by \cite{fordhudson}).
Thus we may assume that $x>199000$. Then it is seen that
$\lambda_{g_{4,3}}(x)\ge 0.4594\mu_{g_{4,3}}(x/59)\ge 0.50456\mu_{g_{3,1}}(x/5)\ge \lambda_{g_{3,1}}(x)$.\\
$\lambda_{g_{4,3}}(x)$ {\tt versus} $\lambda_{g_{4,1}}(x)$.
For $7\le x\le 1.1*10^6$ one directly verifies the inequality (note
that Corollary 2 cannot be used this time). For $x>1.1*10^6$ one
deduces, proceeding as before, that
$\lambda_{g_{4,3}}(x)\ge 0.48508\mu_{g_{4,3}}(x/127)\ge 0.50456\mu_{g_{4,1}}(x/5)\ge \lambda_{g_{4,1}}(x)$.\\

\noindent It remains to establish Theorem \ref{drietje}.\\

\noindent {\it Proof of Theorem} \ref{drietje}. 
We only deal with $N(x;4,3)$ versus $N(x;4,1)$, the other cases
following at once from Theorem \ref{sterker} and (\ref{mlambda}).
Put $\delta(x)=\lambda_{g_{4,3}}(x)-\lambda_{g_{4,1}}(x)$.
By Theorem \ref{sterker} we have $\delta(x)\ge 0$ for
$x\ge 7$. Using this and (\ref{mlambda}) we infer that 
\begin{eqnarray}
N(x;4,3)-N(x;4,1)&=&{\delta(x)\over \log x}
+\int_2^7{\delta(t)dt\over t\log^2 t}
+\int_7^x{\delta(t)dt\over t\log^2 t}\nonumber\\
&\ge& \int_2^7{\delta(t)dt\over t\log^2 t}={\log 5-\log 3\over \log 7 }>0,\nonumber
\end{eqnarray}
for $x\ge 7$. For $x<7$ the result is clearly true. \qed\\

\noindent {\bf Acknowledgement}. T. Dokshitzer was so kind as to redo some
of the computations (carried out in Maple and Quickbasic) in PARI. K. Ford, 
R. Hudson
and M. Rubinstein helpfully provided me with some numerical
data regarding Chebyshev's bias (and with a preprint of \cite{fordhudson}).
I'd like to thank O.
Ramar\'e for pointing out reference \cite{dusart}. He also informed me
that he is working on a paper \cite{ramare0} in which he is able
to roughly half the values given in (\ref{beperking}) (which is not
enough to give an unconditional proof of any of the claims in Theorem 6).
H. te Riele kindly did the computations required to validate
parts a and b of Theorem \ref{cplusafschat}, taking respectively
20 and 5 minutes CPU time on a 300MHz SGI processor.


\begin{thebibliography}{99}
\bibitem{bays} C. Bays and R.H. Hudson, Details of the first region of
integers $x$ with $\pi_{3,2}(x)<\pi_{3,1}(x)$, {\it Math. Comp.}
{\bf 32} (1978), 571-576.
\bibitem{chebyshev} P.L. Chebyshev, Lettre de M. le Professeur
Tch\'ebychev \`a M. Fuss sur un nouveaux th\'eor\`eme relatif aux nombres
premiers contenus dans les formes $4n+1$ ets $4n+3$, {\it Bull. Classe
Phys. Acad. Imp. Sci. St. Petersburg} {\bf 11} (1853), 208.
\bibitem{cohenbook} H. Cohen, {\it Advanced topics in computational number theory}, GTM {\bf 193}, Springer-Verlag, New York, 2000.
\bibitem{cohenpreprint} H. Cohen,
High precision computation of Hardy-Littlewood constants,
draft of a preprint, available at
http://www.math.u-bordeaux/${}^{\sim}$cohen/.
\bibitem{cohen}  H. Cohen and F. Dress, Estimations num\'eriques du reste de la fonction sommatoire relative
aux entiers sans facteur carr\'e,
{\it Publ. Math. Orsay} 88/02 (1988), 73-76.
\bibitem{conway} J.H. Conway and N.J.A. Sloane, Sphere packings, lattices and 
groups, Third edition,
Springer-Verlag, New York, 1999.
\bibitem{cox} D.A. Cox, The arithmetic-geometric mean of Gauss,
{\it Enseign. Math.} {\bf 30} (1984), 275-330.
\bibitem{davenport} H. Davenport, {\it Multiplicative number theory}, 
Third revised edition, Springer-Verlag, New York, 2000.
\bibitem{dilcher} K. Dilcher, Generalized Euler constants for
arithmetical progressions, {\it Math. Comp.} {\bf 59} (1992), 259-282.
\bibitem{dusart} P. Dusart, Autour de la fonction qui
compte le nombre de nombres premiers, PhD thesis, Universit\'e de Limoges, 1998. \bibitem{finch} S. Finch, Mathematical constant web pages,\\
http://www.mathsoft.com/asolve/constant/constant.html
\bibitem{fordhudson} K. Ford and R.H. Hudson, Sign changes
in $\pi_{q,a}(x)-\pi_{q,b}(x)$, submitted for publication.
\bibitem{landau0} E. Landau, \"Uber die zu einem algebraischen Zahlk\"orper
geh\"orige Zetafunktion und die Ausdehnung der Tschebyschefschen Primzahlentheorie aus das Problem der Verteilung der Primideale,
{\it J. reine angew. Math.} {\bf 125} (1903), 64-188.
\bibitem{landau1} E. Landau, \"Uber die Einteilung der positiven ganzen
Zahlen in vier Klassen nach der mindest Anzahl der zu ihrer additiven
Zusammensetzung erforderlichen Quadrate, {\it Arch. der Math.
und Phys.} {\bf 13} (1908), 305-312.
\bibitem{littlewood} J.E. Littlewood, Distribution des nombres premiers,
{\it C.R. Acad. Sci. Paris} {\bf 158} (1914), 1869-1872.
\bibitem{martin} G. Martin, Asymmetries in
the Shanks-R\'enyi Prime Number Race, proceedings of
the Millennial Conference on Number Theory (Urbana, IL), to appear.
\bibitem{blackholes} S.D. Miller and G. Moore, Landau-Siegel zeroes and black hole entropy, {\it Asian J.
Math.} {\bf 4} (2000), 183-211.
\bibitem{moree} P. Moree, Approximation of singular series
and automata, {\it Manuscripta Math.} {\bf 101} (2000), 385-399.
\bibitem{cazaran} P. Moree and J. Cazaran, On a claim of
Ramanujan in his first
letter to Hardy, {\it Exposition.
Math.} {\bf 17} (1999), 289-311.
\bibitem{hexagonal} P. Moree and H.J.J. te Riele, The hexagonal versus
the square lattice, in preparation.
\bibitem{moser} L. Moser and R.A. MacLeod, The error term for the squarefree
integers, {\it Canad. Math. Bull.} {\bf 9} (1966), 303-306.
\bibitem{narkiewicz} W. Narkiewicz, {\it The development of prime number theory.
From Euclid to Hardy and Littlewood}, Springer Monographs in Mathematics,
Springer-Verlag, Berlin, 2000.
\bibitem{nevanlinna} V. Nevanlinna, On constants connected with the prime number theorem for arithmetic progressions, {\it Ann. Acad.
Sci. Fenn. Ser. A I} {\bf 539} (1973), 11 pp..
\bibitem{postnikov} A.G. Postnikov, {\it Introduction to analytic
number theory}, AMS translations of mathematical monographs 68,
AMS, Providence, Rhode Island, 1988.
\bibitem{ramare0} O. Ramar\'e, Sur un th\'eor\`eme de Mertens, in
preparation.
\bibitem{ramare} O. Ramar\'e and R. Rumely, Primes in arithmetic
progressions, {\it Math. Comp.} {\bf 65} (1996), 397-425.
\bibitem{rossers} J.B. Rosser and L. Schoenfeld, Approximate formulas for
some functions of prime numbers, {\it Illinois Journal Math.} {\bf 6}
(1962), 64-94.
\bibitem{chebbie} M. Rubinstein and P. Sarnak, Chebyshev's bias,
{\it Experiment. Math.} {\bf 3} (1994), 173-197.
\bibitem{schmutz} P. Schmutz Schaller, Geometry of Riemann surfaces
based on closed geodesics, {\it Bull. Amer. Math. Soc. (N.S.)}
{\bf 35} (1998), 193-214.
\bibitem{schoenfeld} L. Schoenfeld, Sharper bounds for the Chebyshev
functions $\theta(x)$ and $\psi(x)$. II, {\it Math. Comp.} {\bf 30}
(1976), 337-360.
\bibitem{serre} J.-P. Serre, Divisibilit\'e de certaines fonctions
arithm\'etiques, {\it Enseign. Math.} {\bf 22} (1976), 227-260.
\bibitem{shanks} D. Shanks, The second-order term in the
asymptotic expansion of
$B(x)$, {\it Math. Comp.} {\bf 18} (1964), 75-86.
\bibitem{song} J.M. Song, Sums over nonnegative multiplicative functions
over integers without large prime factors. I, {\it Acta Arith.}
{\bf 97}, 329-351.
\bibitem{turanbook} P. Tur\'an, {\it On a new method of analysis and its applications}, John Wiley and Sons, Inc., New York, 1984.
\bibitem{turan} P. Tur\'an, Collected papers of Paul Tur\'an, Vol. 1-3,
 Ed. P. Erd\"os. Akad\'emiai Kiad\'o, Budapest, 1990.
\bibitem{williams} K.S. Williams, Mertens' theorem for arithmetic
progressions, {\it J. Number Theory} {\bf 6} (1974), 353-359.
\bibitem{wintner} A. Wintner, On the distribution of the remainder
term of the prime number theorem, {\it Amer. J. Math.} {\bf 63}
(1941), 233-248.
\bibitem{wirsing2} E. Wirsing, Das asymptotische Verhalten von
Summen \"uber multiplikative Funktionen, {\it Math. Ann.}
{\bf 143} (1961), 75-102.


\end{thebibliography}
\end{document}